\documentclass[12pt]{article}
\usepackage{amsmath}
\usepackage{amssymb}
\newtheorem{theorem}{Theorem}
\begin{document}
\title{Deformations of Lie algebras of type ${D}_{l}$ and $\overline{{D}_{l}}$ in characteristic 2}

\author{Chebochko N.G.} \footnotetext{The research was supported by the Laboratory of Dynamical Systems and Applications of NRU
HSE, the Ministry on Science and Higher Education of the Russian Federation, grant no. 075-15-
2019-1931, with the exception of the results of Theorems 2 (2). The results of Theorems 2 (2) 
are obtained with the support of RFBR according to research project N 18-01-00900.}
\date{}
\maketitle
\begin{center}
{\it {National Research University Higher School of Economics,\\  Nizhny Novgorod, Russia}}\\
e-mail:
chebochko@mail.ru
\end{center}

\begin{abstract}
 The study of global deformations of Lie algebras is related to the problem of classification of simple Lie algebras over fields of small characteristic. The classification of finite-dimensional simple Lie algebras is complete over algebraically closed fields of characteristic $p>3$ (\cite{Strade1}-\cite{Strade3}). Over the fields of characteristic 2, a large number of examples of Lie algebras are constructed that do not fit into previously known schemes. Description of the deformation of classical Lie algebras first gives new examples of simple Lie algebras, and second allows you to describe known examples as deformations of classical Lie algebras, as is done, for example, in \cite{CKK}. This paper describes the spaces of local deformations of Lie algebras of the type $D_l$ for  $l>3$ and the factor of the algebra at the center $\overline{D}_l$, and describes some global deformations of Lie algebras of these types that are not rigid.
\end{abstract}

Classical deformation theory of associative and Lie algebras
begins with the works of M. Gerstenhaber \cite{gerstenhaber} and
A. Nijenhuis, R. Richardson \cite{nijenhuis} in 1960. They
we studied single parameter deformations and established a connection
between Lie algebra cohomology and infinitesimal cohomology
deformations.

The problem of describing deformations of the Lie algebra is divided into two
tasks:

1) describe the space of local deformations, the elements of which are
are cocycles from the second group of cohomology of the Lie algebra
with coefficients in the conjoined module;

2) solve the problem of the continuity of local deformations. For
the solution to this problem is to find a split into orbits
spaces of local deformations with respect to the action
the group of automorphisms of the lie algebra and test
integrability of representatives from each orbit.

The space of local deformations or the second group of cohomology of the Lie algebra $L$ describes possibilities for constructing deformations of the Lie algebra. The triviality of the second group of cohomology $H^2 (L,L)$ means that there is a neighborhood of the Lie algebra $L$ on the variety (c
zarissky topology) of structures of Lie algebras on
vector space $L$, all points of which are algebras
Whether, isomorphic to $L$.  The Lie algebra in this case is rigid.

A classical Lie algebra over a field of characteristic $p>0$ is a Lie algebra obtained by reduction modulo $p$ of the Chevalley order of complex simple Lie algebra, or its quotient algebra by the centre. Let $L$ be Lie algebra over field $K$, ${{L}_{t}}$  denote the vector space obtained from $L$ by extending the coefficient domain from $K$ to $K((t))$, i.e., $L\otimes K((t))$. Suppose that bilinear function ${{f}_{t}}:{{L}_{t}}\times {{L}_{t}}\to {{L}_{t}}$ has the form

${{f}_{t}}(x,y)=[x,y]+t{{F}_{1}}(x,y)+{{t}^{2}}{{F}_{2}}(x,y)+\cdots$

where ${{F}_{i}}$ -- bilinear function over a field $K$.
If ${{f}_{t}}$ satisfies condition of antisymmetry and Jacobi identity, then Lie algebras ${{L}_{t}}$ with multiplication ${{f}_{t}}$ are family of global deformations of Lie algebra $L$, and ${{F}_{1}}$ is called integrable. Conditions on ${{f}_{t}}$ mean, in particular, that ${{F}_{1}}$ is an element of the group ${{Z}^{2}}(L,L)$. For cohomologous cocycles corresponding Lie algebras are isomorphic. So we can choose representatives of cohomological classes as ${{F}_{1}}$. Second cohomology group with coefficients in adjoint module ${{H}^{2}}(L,L)$ is called space of local deformations. Lie algebra, which does not have non-trivial global deformations, is called rigid.

Classical Lie algebras over fields of characteristic 2 and 3 admit global deformations. Over a field of characteristic 3 there exist non-isomorphic Lie algebras with root system of type ${C}_{2}$. A. I. Kostrikin in \cite{kostrikin} constructed parametric families of non-isomorphic simple Lie algebras over a field of characteristic 3, which are global deformations of Lie algebra of type ${C}_{2}$.  A. S. Dzhumadildaev in \cite{Dzum} has proved that among Lie algebras of type ${A}_{n}$, ${B}_{n}$, ${C}_{n}$, ${D}_{n}$ over a field of characteristic 3 only ${C}_{2}$ admits non-trivial global deformations. A. I. Kostrikin and M. I. Kuznetsov in \cite{kostrikinkuz} fully described global deformations of Lie algebra of type ${C}_{2}$ over a field of characteristic 3. Their approach for study of global deformations, based on the study of orbits of the action of the automorphism group of the algebra on its cohomology space, is used in this paper. A.N. Rudakov \cite{rudakov} has proved that over a field of characteristic $p>3$ all classical Lie algebras are rigid. In \cite{kiril} and \cite{natimi} M. I. Kuznetsov and N.G. Chebochko proposed a new scheme for studying rigidity and proved that over a field of characteristic $p>2$ all classical Lie algebras are rigid except ${C}_{2}$ when $p=3$. Spaces of local deformations of classical Lie algebras with homogeneous root system over a field of characteristic 2 are described in \cite{chebochko}, in particular are described spaces of local deformations of classical Lie algebra of type ${{D}_{l}}$ over a field of characteristic 2.

Let $R$ be root system of type ${{D}_{l}}$, $\{{{\alpha }_{1}},\ldots ,{{\alpha }_{l-2}},{{\alpha }_{l-1}},{{\alpha }_{l}}\}$ -- simple roots, where enumeration is chosen in such way that ${{\alpha }_{l}}$ is connected to ${{\alpha }_{l-2}}$ on Dynkin diagram. By $<\alpha ,\beta >$ we denote Cartan number for roots $\alpha ,\beta $. By ${{H}_{i}}={{H}_{{{\alpha }_{i}}}}$ and ${{E}_{\alpha }}$, $\alpha \in R$ we denote vectors of Chevalley basis in Lie algebra.

\section{Deformations of Lie algebra of type ${\bar{D}}_{l}$ for odd $l>3$.}

For odd $l>3$ Lie algebra of type   over a field $K$ of characteristic 2 has one-dimensional centre ${{H}_{l}}+{{H}_{l-1}}$. Let L be Lie algebra of type ${\bar{D}}_{l}$ --  quotient algebra of   by the centre.
According to \cite{frohardt}, the automorphism group of Lie algebra $L$ of type ${\bar{D}}_{l}$ is isomorphic to Chevalley group of type ${B}_{l}$. Let $V$ be vector space over $K$ of dimension $2l$ with a symplectic form $(\ ,\,\ )$. The automorphism group of $L$ contains group $G=Sp(2l)$ --  the symplectic group, associated with the form$(\ ,\,\ )$. We choose a symplectic basis $\{{{e}_{1}},\,\ \ldots ,\,\ {{e}_{l}},\,\ {{e}_{-l}},\,\ \ldots ,\,\ {{e}_{-1}}\}$  in  $V$, consisting of eigenvectors
under the action of a maximal torus $T$ of the group $G$. Let  ${e}_{i}$  have weight ${{\varepsilon }_{i}}$,  ${e}_{-i}$  have weight $-{{\varepsilon }_{i}}$,  $i=1,\ldots ,\,l$. The quotient algebra of
the exterior algebra $\Lambda^{2}V$ by the ideal $I=<{{e}_{1}}{{e}_{-1}}+\cdots +{{e}_{l}}{{e}_{-l}}>$ is isomorphic to the Lie algebra $L$ (this is true only for odd $l$).
The basic complex of cohomologies can be decomposed into the direct sum of weight subcomplexes using the natural action of a maximal torus of the Chevalley group G(L). The corresponding cohomology groups are weight subspaces in the cohomology group of the basic complex. We denote by $C_{\mu }^{n}(L,L)$, $Z_{\mu }^{n}(L,L)$, $B_{\mu }^{n}(L,L)$, $H_{\mu }^{n}(L,L)$ the weight subspaces of cochains, cocycles, coboundaries and cohomologies of weight $\mu$, respectively.
In \cite{chebochko} it is proved that $H_{0}^{2}(L,L)=0$, ${H}^{2}(L,\,L)$  is the direct sum of non-zero subspaces  $H_{\mu }^{2}(L,\,L)$ where $\mu $ has the form $\gamma +\delta $ for $\gamma ,\ \delta \in R$ satisfying condition $<\gamma ,\delta >=0$. Any weight $\mu $ is conjugate under the action of the Weyl group with ${{\alpha }_{l}}+{{\alpha }_{l-1}}$. Total number of such weights is $2l$: $\{\pm ({{\alpha }_{l}}-{{\alpha }_{l-1}}),\ \pm ({{\alpha }_{l}}+{{\alpha }_{l-1}}),\pm ({{\alpha }_{l}}+{{\alpha }_{l-1}}+2{{\alpha }_{l-2}}),\ \ldots ,\ \pm ({{\alpha }_{l}}+{{\alpha }_{l-1}}+2{{\alpha }_{l-2}}+\cdots +2{{\alpha }_{1}})\}$ and $\dim{{H}^{2}}(L,\,L)=2l$. Also natural description of  cohomology group ${H}^{2}(L,\,L)$ as a module over the group  $Sp(2l)\subset $$Aut(L)$ is presented in \cite{chebochko}. It is proved that for Lie algebra of type ${\bar{D}}_{l}$  ( $l>3$ - odd)  ${H}^{2}(L,\,L)$  as a module over  $Sp(2l)$  is, up to the Frobenius morphism, isomorphic to $V$, where $V$ is the standard module over $Sp(2l)$. Below we describe in details this implementation of ${H}^{2}(L,\,L)$.
Lets construct a homomorphism of the module $V$ into the module ${H}^{2}(L,\,L)$.
First we define $\Phi \ \ :\ \ V\to {{C}^{2}}(L,\,L)$. Let $v\in V$,  ${w}_{1}{{w}_{2}},\,\ {{w}_{3}}{{w}_{4}}\in \Lambda^{2}V/I\cong L$. Set $\Phi (v)=\varphi $, where

$\varphi ({{w}_{1}}{{w}_{2}},{{w}_{3}}{{w}_{4}})=(v,\,{{w}_{1}})(v,\,{{w}_{3}}){{w}_{2}}{{w}_{4}}+(v,\,{{w}_{2}})(v,\,{{w}_{3}}){{w}_{1}}{{w}_{4}}+ \\
 +(v,\,{{w}_{1}})(v,\,{{w}_{4}}){{w}_{2}}{{w}_{3}}+(v,\,{{w}_{2}})(v,\,{{w}_{4}}){{w}_{1}}{{w}_{3}}+(v,\,{{w}_{1}})({{w}_{3}},\,{{w}_{4}})v{{w}_{2}}+ \\
$
$+(v,\,{{w}_{2}})({{w}_{3}},\,{{w}_{4}})v{{w}_{1}}+(v,\,{{w}_{3}})({{w}_{1}},\,{{w}_{2}})v{{w}_{4}}+(v,\,{{w}_{4}})({{w}_{1}},\,{{w}_{2}})v{{w}_{3}}.$

Using the Poisson brackets: $\{{{v}_{1}}{{v}_{2}},\,v\}=({{v}_{1}},\,v){{v}_{2}}+({{v}_{2}},\,v){v}_{1}$, $\Phi (v)$ can be defined as follows:

$\Phi (v)({{w}_{1}}{{w}_{2}},{{w}_{3}}{{w}_{4}})=\{\,{{w}_{1}}{{w}_{2}},v\}\{\,{{w}_{3}}{{w}_{4}},v\}+v\{v,(\,{{w}_{3}},{{w}_{4}})\,
{{w}_{1}}{{w}_{2}}+\\(\,{{w}_{1}},{{w}_{2}})\,{{w}_{3}}{{w}_{4}}\}$.

$\varphi $ is correctly defined on elements of $\Lambda^{2}V$, also $\varphi $ is linear and
skew-symmetric, and, consequently, is an element of ${C}^{2}(L,\,L)$. Also $\Phi $ commutes with the action of $Sp(2l)$:

$\Phi (gv)({{w}_{1}}{{w}_{2}},\,{{w}_{3}}{{w}_{4}})=g(\Phi (v)({{g}^{-1}}{{w}_{1}}{{g}^{-1}}{{w}_{2}},\,{{g}^{-1}}{{w}_{3}}{{g}^{-1}}{{w}_{4}}))$

for any $g\in Sp(2l)$.
Using the Poisson brackets in  $\Lambda^{2}V$:

$[{{v}_{1}}{{v}_{2}},\,{{v}_{3}}{{v}_{4}}]=\{{{v}_{1}}{{v}_{2}},\,{{v}_{3}}{{v}_{4}}\}=({{v}_{1}},\,{{v}_{3}}){{v}_{2}}{{v}_{4}}+({{v}_{1}},\,{{v}_{4}}){{v}_{2}}{{v}_{3}}+({{v}_{2}},\,{{v}_{3}}){{v}_{1}}{{v}_{4}}+({{v}_{2}},\,{{v}_{4}}){{v}_{1}}{v}_{3}$,                           calculation of $d\varphi $ on arbitrary elements ${w}_{1}{w}_{2}$, ${w}_{3}{w}_{4}$, ${w}_{5}{w}_{6}$ shows that $\varphi \in {{Z}^{2}}(L,L)$.
Mapping $\Phi $ induces mapping of $V$ into cohomology group ${{H}^{2}}(L,L)$, which we will also denote as $\Phi $. So, $\Phi ({{e}_{\pm 1}})$,  $\ldots $ ,$\Phi ({{e}_{\pm l}})$  define  $2l$  non-trivial cocycles of different weights. Since $\Phi (kv)={{k}^{2}}\Phi (v)$ for any $k\in K$, the action of symplectic group on ${{H}^{2}}(L,L)$ is the action on $V$ up to the Frobenius morphism.

The necessary condition of integrability of arbitrary cocycle $\psi $ to global deformation is the triviality of cocycle $\psi \cup \psi $  from  ${Z}^{3}(L,\,L)$, where

 $\psi \cup \psi (x,\,y,\,z)=\psi (\psi (x,\,y),\,z)+\psi (\psi (y,\,z),\,x)+\psi (\psi (z,\,x),\,y)$

for any $x,y,z$ from $L$.
In Lie algebra of type ${\bar{D}}_{l}$ ($l>3$  - odd) exists cocycle  $\psi $, such that  $\psi \cup \psi \not=0$  and  $\psi \cup \psi $  has weight  $\mu $,  for which  $C_{\mu }^{2}(L,\,L)=0$. Lets prove it.
Since $l>3$  is odd, in basis $\{{{e}_{1}},\,\ \ldots ,\,\ {{e}_{l}},\,\ {{e}_{-l}},\,\ \ldots ,\,\ {{e}_{-1}}\}$ exist vectors ${{e}_{4}},{{e}_{5}},{{e}_{-5}},{{e}_{-4}}$. Consider the cocycle $\psi =\Phi ({{e}_{4}})$. We have

$   \psi \cup \psi ({{e}_{-4}}{{e}_{-5}},\,{{e}_{-4}}{{e}_{5}},\,{{e}_{3}}{{e}_{-4}})=\psi (\psi ({{e}_{-4}}{{e}_{-5}},\,{{e}_{-4}}{{e}_{5}}),\,{{e}_{3}}{{e}_{-4}})+\\ \psi (\psi ({{e}_{-4}}{{e}_{5}},\,{{e}_{3}}{{e}_{-4}}),{{e}_{-4}}{{e}_{-5}})+
 \psi (\psi ({{e}_{3}}{{e}_{-4}},{{e}_{-4}}{{e}_{-5}}),\,{{e}_{-4}}{{e}_{5}}\,)=\\\psi (\psi ({{e}_{-4}}{{e}_{-5}},\,{{e}_{-4}}{{e}_{5}}),\,{{e}_{3}}{{e}_{-4}})=\psi ({{e}_{5}}{{e}_{-5}},\,{{e}_{3}}{{e}_{-4}})={{e}_{3}}{{e}_{4}}. \\
$

It is seen that $\psi \cup \psi $ is nonzero cocycle from $Z_{4{{\varepsilon }_{4}}}^{3}(L,L)$. Since space $C_{4{{\varepsilon }_{4}}}^{2}(L,L)$ is zero, $B_{4{{\varepsilon }_{4}}}^{3}(L,L)=0$. Therefore, $\psi \cup \psi $ is not a coboundary.
So, $\psi $ is nonintegrable cocycle.
The group $G$ acts transitively on $V$, so  cocycle from ${H}^{2}(L,\,L)$  is not integrable and Lie algebras of type ${\bar{D}}_{l}$  ( $l>3$ - odd) are rigid.

\section{Deformations of Lie algebra of type ${D}_{l}$ for even $l\ge 4$.}

Let L be Lie algebra of type ${D}_{l}$ for even $l\ge 4$. Then $L$ has two-dimensional centre: $<{{H}_{l}}+{{H}_{l-1}},\ {{H}_{l-1}}+{{H}_{l-3}}+\cdots +{{H}_{3}}+{{H}_{1}}>$. The space ${H}^{2}(L,\,L)$  is the direct sum of non-zero subspaces $H_{\mu }^{2}(L,\,L)$, where  $\mu $  has the form $\gamma +\delta $ for $\gamma ,\ \delta \in R$, satisfying condition $<\gamma ,\delta >=0$. For $l>4$ any of such weights $\mu $ is conjugate under the action of the Weyl group with ${{\alpha }_{l}}+{{\alpha }_{l-1}}$. Total number of such weights is $2l$: $\{\pm ({{\alpha }_{l}}-{{\alpha }_{l-1}}),\ \pm ({{\alpha }_{l}}+{{\alpha }_{l-1}}),\pm ({{\alpha }_{l}}+{{\alpha }_{l-1}}+2{{\alpha }_{l-2}}),\ \ldots ,\ \pm ({{\alpha }_{l}}+{{\alpha }_{l-1}}+2{{\alpha }_{l-2}}+\cdots +2{{\alpha }_{1}})\}$. Corresponding weight subspaces of cohomologies are one-dimensional and $\dim{{H}^{2}}(L,\,L)=2l$.
For $l=4$ there are 24 weights and they are conjugate with ${{\alpha }_{1}}+{{\alpha }_{3}}$, ${{\alpha }_{1}}+{{\alpha }_{4}}$, ${{\alpha }_{3}}+{{\alpha }_{4}}$, $\dim{{H}^{2}}(L,\,L)=24$.
We will use isomorphism ${{C}^{n}}(L,L)\cong \underbrace{{{L}^{*}}\wedge \ldots \wedge {{L}^{*}}}_{n}\otimes L$.
Since all weights are conjugate (in case of ${D}_{4}$ conjugate with ${{\alpha }_{1}}+{{\alpha }_{3}}$, ${{\alpha }_{1}}+{{\alpha }_{4}}$, ${{\alpha }_{3}}+{{\alpha }_{4}}$), it is sufficient to describe space $H_{{{\alpha }_{l}}+{{\alpha }_{l-1}}}^{2}(L,\,L)$.
Cocycle

$\psi =\sum\limits_{\gamma +\delta ={{\alpha }_{l}}+{{\alpha }_{l-1}}}{E_{-\gamma }^{*}}\wedge E_{-\delta }^{*}\otimes z$,

where $z={{H}_{l-1}}+{{H}_{l-3}}+\cdots +{{H}_{3}}+{{H}_{1}}$, generates subspace $H_{{{\alpha }_{l}}+{{\alpha }_{l-1}}}^{2}(L,\,L)$.
Condition integrability for cocycle $\psi $ has the form:

$\psi \cup \psi (x,\,y,\,z)=0$                                                                                                                  for any $x,y,z$ from $L$.
Therefore, cocycle $\psi $ is integrable. Mapping
${{f}_{t}}(x,y)=[x,y]+t\psi (x,y)$
gives global deformation of Lie algebra $L$.
Since for any cocycle $\psi \in {{H}^{2}}(L,\,L)$ $\psi (x,\,y)$ is contained in the centre of $L$ for any $x,y$ from $L$ and $\psi(x,y)=0$ if $x,y \in Z(L)$, condition $\psi \cup \psi (x,\,y,\,z)=0$ holds true for all elements from ${{H}^{2}}(L,\,L)$. So any cocycle ${{H}^{2}}(L,\,L)$ for Lie algebra of type ${D}_{l}$ for even $l$ are integrable.

{\bf Conclusions.}

Main results of this paper are formulated in the following theorems.
\begin{theorem}
Let $L$ be Lie algebra over an algebraically closed field of characteristic 2.

(1)	If $L$ has type ${D}_{l}$ for for odd $l>3$ and ${\bar{D}}_{l}$ for even $l\neq 6$  then $L$ have a trivial space of local deformations and are rigid.

(2) If $L$ has type ${D}_{l}$ for even $l\neq 4$ and ${\bar{D}}_{l}$ for odd $l$  then the space of local deformations has dimension $2l$.

(3)	The space of local deformations of the Lie algebra $ \overline{D}_6$ has dimension 64. The space of local deformations of the Lie algebra ${D}_4$ has dimension 24.

\end{theorem}

\begin{theorem}
Let $L$ be Lie algebra over an algebraically closed field of characteristic 2.

(1)	If $L$ has type ${\bar{D}}_{l}$ for odd $l>3$  then all cocycles of the group $H^2 (L,L)$ lie in the same orbit relative to the action of the automorphism group and they are not integrable. $L$ is rigid Lie algebra.

(2)	If $L$ has type ${D}_{l}$ for even $l\ge 4$, then any cocycle  $\psi\in {{H}^{2}}(L,\,L)$ is integrable and defines a global deformation of the form $f_t=[~,~]+t\psi$.

\end{theorem}

{ \bf Acknowledgments}

The author is grateful to M. I. Kuznetsov for attention to this work and useful remarks.

\end{document}